\documentclass[12pt]{article}
\usepackage{rotating}
\usepackage{amsfonts}
\usepackage[round]{natbib}

\title{A note on Horvitz-Thompson estimators for rare subgroup analysis in the presence of interference}

\author{Erin E Gabriel\\
Department of Medical Epidemiology and Biostatistics,\\ Karolinska Institutet, Stockholm, Sweden\\
e-mail: erin.gabriel@ki.se}

\begin{document}

\maketitle

\noindent ABSTRACT\\
When there is interference, a subject's outcome depends on the treatment of others and treatment effects may take on several different forms. This situation arises often, particularly in vaccine evaluation. In settings where interference is likely, two-stage cluster randomized trials have been suggested as a means of estimating some of the causal contrast of interest. Working in the finite population setting to investigate rare and unplanned subgroup analyses using some of the estimators that have been suggested in the literature, include Horvitz-Thompson, H\'ajek, and what might be called the natural extension of the marginal estimators suggested in Hudgens and Halloran 2008. I define the estimands of interest conditional on individual, group and both individual and group baseline variables, giving unbiased Horvitz-Thompson style estimates for each. I also provide variance estimators for several estimators. I show that the Horvitz-Thompson (HT) type estimators are always unbiased provided at least one subject within the group or population, whatever the level of interest for the estimator, is in the subgroup of interest. This is not true of the ``natural" or the H\'ajek style estimators, which will often be undefined for rare subgroups. \\

\noindent \textbf{keywords}: Finite population; Horvitz-Thompson estimators; Interference.

\section{Introduction}
Interference is present in settings where the subject's outcomes are not independent of other subject's treatment. Interference changes both the types of causal effects of treatment and their estimation. Much of the causal inference literature assumes that units of interest are independent. In the presence of interference, causal inference becomes much more difficult. However, interest has increased in the last decade with regards to causal inference in the presence of interference. There have been several papers which discuss causal inference with interference starting with the foundational work of \citet{hudgens08} which was extended by \citet{tchetgen10}, and since then there has been an explosion of new works \citet{liu14,vanderweele13, liu2016inverse, aronow2017estimating, savje2017average}. Many papers focus on marginal effects and estimands, \citet{vanderweele12} and \citet{vanderweele13} deal with conditional estimands as does \citet{halloran12}. 

I am interested in randomization inference within subgroups defined by rare baseline variables where the randomization has ignored the baseline variable(s) in question and may therefore not contain any members of the subgroup in one of the randomized arms. In this setting, what one might call the ``natural" estimators or the H\'ajek estimators will be either undefined or biased. If one could randomize stratified by the desired conditioning variable, all \citet{hudgens08} and \citet{tchetgen10} theorems and propositions would apply directly. Although this may seem like a minor point, pre-specification of intended analysis is often required in randomized clinical trials, thus one cannot change the analysis after looking at the unblinded data. 

I show that the Horvitz-Thompson estimators \citet{Horvitz52} are unbiased and defined, provided there are any members of the subgroup in the population of interest, regardless of the realized randomization. I consider two-stage randomization to a fixed number of clusters at the first stage, and then a fixed number of subjects within the cluster at the second stage \citet{hudgens08}. For simplicity, I make the same assumptions as \citet{hudgens08}. However, as has now been shown now in many works, there are estimators that are unbiased for interesting estimands under reduced assumptions \citet{savje2017average, aronow2017estimating}. 

\section{Notation} \label{NO}
Following the notation of \citet{hudgens08}, suppose there are $J>1$ clusters of individuals. For $j \in [1, \ldots, J]$, let $n_j$ denote the number of subjects in the $j$th cluster, indexed by $i$. Let the treatment assignments for the individuals in cluster $j$ be denoted by the vector $\mathbf{Z}_{j} \equiv (Z_{j1}, \ldots, Z_{jn_{j}})$, with the sub-vector excluding the treatment assignment of subject $i$ written as $\mathbf{Z}_{j(i)}$. For simplicity of notation and illustration, the treatment assignment for subject $i$ in cluster $j$ is either 0 or 1 with the realization being denoted by $z_{ij}$. The realization of the cluster level treatment assignment $\mathbf{z}_{j}$ can be any of $2^{n_{j}}$ possible values; realizations excluding the assignment for subject $i$ will be written as $\mathbf{z}_{j(i)}$. 

Let $d_j$ be some baseline variable(s) observed at the cluster level prior to randomization and $w_{ij}$ be some, not necessarily related, baseline variable(s) at the individual level. Let these variables have an arbitrary domain, and let $B$ or $b$ denote some subset of that domain at the cluster or individual level, respectively. I will consider two coverage proportions for the same intervention, denoted by $\alpha$ and $\gamma$. I define the vector of cluster assignments to these strategies as $\textbf{Q}$, where $Q_{j} =1$ if cluster $j$ is assigned to strategy $\alpha$ and is $0$ otherwise. 

\noindent
Let $m_{j,b,z} = \sum^{n_j}_{i=1} \mathbb{I}_{[z_{ij}=z]}\mathbb{I}_{[w_{ij} \in b]}$ with $z \in \{0,1\}$ and $M_{j,b}=\sum_z m_{j,b,z}$.\\
Let the number of clusters with $M_{j,b}>0$ be denoted by $M_{b}=\sum_j\mathbb{I}_{[M_{j,b}>0]}$.\\
Let $m_{B,q}=\sum^{J}_{j=1} \mathbb{I}_{Q_j=q}\mathbb{I}_{[d_j \in B]}$ for $q \in \{\alpha,\gamma\}$ and $M_{B} =  \sum_{q} m_{B,q}$. \\
Let the number of clusters that have both $M_{j,b}>0$ and $d_j \in B$ 
be denoted by $M_{B,b} = \sum_j  \mathbb{I}_{[M_{j,b}>0]} \mathbb{I}_{[d_j \in B]}$.\\
Let the set of groups with $M_{j,b}>0$, be denoted as $\mathbb{J}^b$, and the set of groups with $d_j \in B$ as $\mathbb{J}^B$.

Let $Y_{ij}(\mathbf{Z}_{j(i)}=\mathbf{z}_{j(i)}, Z_{ij}=z)$ be the potential clinical outcome of subject $i$ in cluster $j$ given the cluster level treatment assignment was $\mathbf{z}_{j(i)}$ and under intervention $z$ for the subject. This allows the $i$th subject's outcome to be influenced by both $Z_{ij}$ and the treatment of other subjects in the same cluster, $\mathbf{Z}_{j(i)}$. Let the realized outcome for subjects $i$ in cluster $j$ be denoted as $Y_{ij}$. 

 I will assume two-part randomization, under which I first randomize clusters to $\alpha$ or $\gamma$ then randomize subjects within each cluster to match the given strategy. The randomization strategy is assumed to be mixed, assigning a set number of clusters and then a set number of subjects within each cluster to treatment as in \citet{hudgens08}. Let the set of all possible randomizations for a cluster of size $n$ that satisfy strategy $\alpha$ be denote as $R^{n}_{\alpha}$ and the subset of these randomizations for which subject $i$ in this cluster is assigned $z$ be denoted by $R^{n-1}_{z;\alpha}$. It should be noted that this notation assumes that there is no interference between clusters. 

Let $Y_{ij}(\mathbf{Z}_{j}=\mathbf{z}_{j})$ denote the potential outcome of subject $i$ in cluster $j$ if the  randomization within the cluster was realized to be $\mathbf{z}_{j}$. I define the individual average clinical outcome under $Z_{ij}=z$ for strategy $\alpha$ as

$$\overline{Y}_{ij}(z; \alpha) = \sum_{\nu \in R^{n_j-1}_{z;\alpha}} Y_{ij}(\mathbf{Z}_{j(i)} = \nu, Z_{ij}=z)   P_{\alpha}(\mathbf{Z}_{j(i)} = \nu\mid  Z_{ij}=z).$$

The cluster average outcome is then given by  $\overline{Y}_{j}(z; \alpha) \equiv 1/n_j \sum^{n_j}_{i=1} \overline{Y}_{ij}(z; \alpha)$. I can also define the marginal potential outcomes, marginalizing over $Z_{ij}$ within a cluster. For subject $i$ in cluster $j$ under strategy $\alpha$ let $\overline{Y}_{ij}(\alpha)$ be the individual average marginal clinical outcome defined by $$\overline{Y}_{ij}(\alpha)= \sum_{\nu \in R^{n_j}_{\alpha}} Y_{ij}(\mathbf{Z}_{j} = \nu)   P_{\alpha}(\mathbf{Z}_{j} = \nu).$$ I can now define the conditional estimands of interest. Based on these one can define cluster and population level summaries of these potential outcomes as in \citet{hudgens08}, as well as contrasts of interest for defining causal effects. 

The group average potential outcomes are $\overline{Y}_{j}(z;\alpha) \equiv 1/n_j\sum^{n_j}_{i=1}\overline{Y}_{ij}(z;\alpha)$ and the population average potential outcomes are
$\overline{Y}(z;\alpha) \equiv \frac{1}{J} \sum^{J}_{j=1} \overline{Y}_{j}(z;\alpha)$
and $\overline{Y}(\alpha) \equiv 1/J \sum^{J}_{j=1} \overline{Y}_{j}(\alpha).$

I can then define contrasts of these potential outcomes to define causal effects. The direct effect, as given in \citet{hudgens08} in a cluster $j$ is given by $DE_j(\alpha) \equiv \overline{Y}_{j}(1;\alpha) -\overline{Y}_{j}(0;\alpha)$
and the direct effect at the population level is defined as $DE(\alpha) \equiv \overline{Y}(1;\alpha) -\overline{Y}(0;\alpha).$

The indirect effect at the population level comparing $\alpha$ and $\gamma$ is defined as
$IE(\alpha, \gamma) \equiv \overline{Y}(0;\alpha) -\overline{Y}(0;\gamma)$ and the population direct effect plus the indirect effect is the total effect, $TE(\alpha, \gamma) \equiv \overline{Y}(1;\alpha) -\overline{Y}(0;\gamma).$ The population overall effect comparing $\alpha$ and $\gamma$ is defined as $OE(\alpha, \gamma) \equiv \overline{Y}(\alpha) -\overline{Y}(\gamma).$ Other definitions of the direct effects, as well as decomposition have been considered \citet{vanderweele11}.

\subsection{Conditional estimands} \label{EST}

I now consider baseline variable conditional versions of the above estimands that are conditional in three ways, conditional on individual level baseline covariates, conditional on cluster level covairates and conditional on both. These are the same estimands as those considered in \citet{hudgens08}, but within a subgroup defined by the baseline variable. Consider that individual $i$ in group $j$ has $w_{ij} \in b$ the individual average outcome under $Z_{ij}=z$ is $\overline{Y}_{ij}(z;\alpha\mid b) \equiv \overline{Y}_{ij}(z;\alpha)$ and is zero otherwise. Similarly, $\overline{Y}_{j}(z;\alpha\mid B) \equiv \overline{Y}_{j}(z;\alpha)$ if cluster $j$ has cluster level baseline variable $d_j \in B$ and is zero otherwise. Finally, the cross-conditional estimand $\overline{Y}_{j}(z;\alpha\mid B,b)$ is equal to  $\overline{Y}_{j}(z;\alpha\mid b)$ if cluster $j$ has cluster level baseline variable $d_j \in B$ and is zero otherwise. Thus, individual and cluster level conditional estimands are special cases of the cross-conditional estimands, when all subjects or all clusters are within the range of interest for the individual or cluster level baseline variable. Hence, this is how they are displayed in Table \ref{estimands1}. I discuss estimators for each type of conditioning, and properties of them, separately for clarity. 

\begin{sidewaystable}[]
    \centering
     \caption{Conditional Estimand Definitions \label{estimands1}}
            \begin{tabular}{l|ccc}
    Notation & \multicolumn{3}{c}{Definitions for given conditioning} \\
    &$\mid B,b$&All $d_i \in B$ $(\mid B,b \equiv \mid b)$ & All $w_{ij} \in b$ $(\mid B,b \equiv \mid B)$ \\
        $\overline{Y}_{ij}(z;\alpha\mid B,b)\equiv$ & $\overline{Y}_{ij}(z;\alpha)\mathbb{I}_{[d_j \in B]}\mathbb{I}_{[w_{ij} \in b]}$  &$\overline{Y}_{ij}(z;\alpha)\mathbb{I}_{[w_{ij} \in b]}$ & $\overline{Y}_{ij}(z;\alpha)\mathbb{I}_{[d_j \in B]}$ \\
        
        $\overline{Y}_{ij}(\alpha \mid  B,b)\equiv$ & $ \overline{Y}_{ij}(\alpha)\mathbb{I}_{[d_j \in B]}\mathbb{I}_{[w_{ij} \in b]}$&  $\overline{Y}_{ij}(\alpha)\mathbb{I}_{[w_{ij} \in b]}$ & $\overline{Y}_{ij}(\alpha)\mathbb{I}_{[d_j \in B]}$ \\
        
        $\overline{Y}_{j}(z;\alpha\mid B,b)\equiv$ & $\mathbb{I}_{[d_j \in B]}\frac{1}{M_{j,b}}\sum^{n_j}_{i=1}\overline{Y}_{ij}(z;\alpha\mid B,b)$&  $\frac{1}{M_{j,b}}\sum^{n_j}_{i=1}\overline{Y}_{ij}(z;\alpha\mid b)$ & $\frac{1}{n_j}\sum^{n_j}_{i=1} \overline{Y}_{ij}(z;\alpha\mid B)$\\
        
        $\overline{Y}_{j}(\alpha\mid B,b)\equiv$ & $\mathbb{I}_{[d_j \in B]}\frac{1}{M_{j,b}}\sum^{n_j}_{i=1}\overline{Y}_{ij}(\alpha\mid B,b)$ & $\frac{1}{M_{j,b}}\sum^{n_j}_{i=1}\overline{Y}_{ij}(\alpha\mid b)$ & $\frac{1}{n_j}\sum^{n_j}_{i=1} \overline{Y}_{ij}(\alpha\mid B)$\\
        
        $DE_j(\alpha\mid B,b)\equiv$ &$\overline{Y}_{j}(1;\alpha\mid B,b) -\overline{Y}_{j}(0;\alpha\mid B,b)$ &  $\overline{Y}_{j}(1;\alpha\mid b) -\overline{Y}_{j}(0;\alpha\mid b)$& $\overline{Y}_{j}(1;\alpha\mid B) -\overline{Y}_{j}(0;\alpha\mid B)$\\
        
        $\overline{Y}(z;\alpha\mid B,b)\equiv$& $\frac{1}{M_{B,b}} \sum^{J}_{j=1} \overline{Y}_j(z;\alpha\mid B,b)$ &  $\frac{1}{M_b} \sum^{J}_{j=1} \overline{Y}_{j}(z;\alpha\mid b)$& $\frac{1}{M_{B}} \sum^{J}_{j=1} \overline{Y}_j(z;\alpha\mid B)$\\
        
        $\overline{Y}(\alpha\mid B,b)\equiv$ & $\frac{1}{M_{B,b}} \sum_J \overline{Y}_j(\alpha\mid B,b)$ &  $\frac{1}{M_b} \sum^{J}_{j=1} \overline{Y}_{j}(\alpha\mid b)$& $\frac{1}{M_{B}} \sum^{J}_{j=1} \overline{Y}_j(\alpha\mid B)$\\
        
        $DE(\alpha \mid b)\equiv$  &$\overline{Y}(1;\alpha\mid B,b) -\overline{Y}(0;\alpha\mid B,b)$ &  $\overline{Y}(1;\alpha\mid b) -\overline{Y}(0;\alpha\mid b)$& $\overline{Y}(1;\alpha\mid B) -\overline{Y}(0;\alpha\mid B)$\\
        
        $TE(\alpha, \gamma\mid B,b)\equiv$& $\overline{Y}(1;\alpha\mid B,b) -\overline{Y}(0;\gamma\mid B,b)$ & $\overline{Y}(1;\alpha\mid b) -\overline{Y}(0;\gamma\mid b)$&  $\overline{Y}(1;\alpha\mid B) -\overline{Y}(0;\gamma\mid B)$\\
        
        $IE(\alpha, \gamma\mid B,b)\equiv$& $\overline{Y}(0;\alpha\mid B,b) -\overline{Y}(0;\gamma\mid B,b)$ & $\overline{Y}(0;\alpha\mid b) -\overline{Y}(0;\gamma\mid b)$&$\overline{Y}(0;\alpha\mid B) -\overline{Y}(0;\gamma\mid B)$\\
        
        $OE(\alpha, \gamma\mid B,b)\equiv$& $\overline{Y}(\alpha\mid B,b) -\overline{Y}(\gamma\mid B,b)$  & $\overline{Y}(\alpha\mid b) -\overline{Y}(\gamma\mid b)$& $\overline{Y}(\alpha\mid B) -\overline{Y}(\gamma\mid B)$\\
    \end{tabular}
    
  \end{sidewaystable}

\section{Assumptions} \label{ASS}
I first need to make assumptions to link the observable data to the desired counterfactuals and estimators.
\begin{itemize}
\item Assumption (a): Consistency,
\item Assumption (b): No interference between clusters,
\item Assumption (c): All assignment strategies are mixed as defined in \citet{hudgens08} and \citet{sobel06}
\item Assumption (d): Stratified interference as defined in \citet{hudgens08}. 
\end{itemize}

Assumption (a) ``consistency" simply means that if $Z_{ij}=z$ and $Q_j =1$ then $Y_{ij}=Y_{ij}(z;\alpha)$, \citep{vanderweele09}. Assumption (b) is implicitly made in the notation and is explained in detail in Section \ref{NO}.  Under Assumption (c), all clusters have the same probability of being assigned a given strategy; $Pr_{\alpha} \equiv P(\mathbf{Z}_j \in\alpha)$ is the probability that a cluster will receive coverage $\alpha$. Under the mixed strategy this is equal to $K/J$, as a fixed number of clusters, $K$ will be randomized to $\alpha$. As well, within a cluster, each individual has the same probability of being assigned to treatment given the randomized coverage which is denoted as $P^{j}_{\alpha} \equiv Pr_{\alpha}(Z_{ij}=1)$, which under a mixed strategy is equal to $k_j/n_j$, as fixed number of subjects, $k_j$ will be randomized to $z_{ij}=1$. Assumption (d), Stratified interference, is an assumption outlined in both \citet{hudgens08}  and \citet{tchetgen10}. It states that only a subject's treatment assignment and the total proportion of people assigned to $z_{ij}=1$ within their cluster impacts their outcome. Therefore, all possible counterfactuals $Y_{ij}(\mathbf{z}_{j(i)}; z_{ij}=z)$ will have the same value for all $\mathbf{z}_{j(i)} \in R^{n_j-1}_{z;\alpha}$.

\section{ESTIMATION AND INFERENCE} \label{thegoods}
I will consider finite population inference following \citet{hudgens08}. Throughout, I mean the expected value $\mbox{E}$ to be with respect to the randomization distribution had each subject and cluster been randomized in each possible pattern, for the fixed and observed, due to consistency, potential outcome. 

Within groups of clusters assigned to $\alpha$, let the cross-conditional outcome estimators be defined by: 
$$\widehat{\overline{Y}}_j(1;\alpha\mid B,b) =\mathbb{I}_{[d_{j} \in B]} \frac{1}{\sum_i \mathbb{I}_{[d_{j} \in b]}}\sum^{n_j}_{i} \frac{Y_{ij}\mathbb{I}_{[Z_{ij}=1]}\mathbb{I}_{[w_{ij} \in b]}}{\mbox{P}^{j}_{\alpha}}$$
with, $\widehat{\overline{Y}}_j(0;\alpha\mid B,b)$ similarly defined, and 
$$\widehat{\overline{Y}}_j(\alpha\mid B,b) = \mathbb{I}_{[d_j \in B]} \frac{1}{\sum_i \mathbb{I}_{[w_{ij} \in b]}}\sum^{n_j}_{i} \frac{Y_{ij} \mathbb{I}_{[w_{ij} \in b]}}{\mbox{P}^{j}_{\alpha}}.$$
This makes the the estimator of $DE_j(\alpha\mid B,b)$,\\ 
$\widehat{DE}_j(\alpha\mid B,b) = \widehat{\overline{Y}}_j(1; \alpha\mid B,b) - \widehat{\overline{Y}}_j(0; \alpha\mid B,b)$.\\

\noindent At the population level, then estimators and then given by: 
$$\widehat{\overline{Y}}(z;\alpha\mid B,b) =\frac{1}{ \sum_j \mathbb{I}_{[d_j \in B]} \mathbb{I}_{[\sum w_{ij} \in b]>0} }\frac{\sum \widehat{\overline{Y}}_j(z;\alpha\mid B,b)\mathbb{I}_{[Q_j=1]}}{\mbox{Pr}_{\alpha}},$$
$$\widehat{\overline{Y}}(\alpha\mid B,b) = \frac{1}{ \sum_j \mathbb{I}_{[d_j \in B]} \mathbb{I}_{[\sum w_{ij} \in b]>0}}\frac{\sum \widehat{\overline{Y}}_j(\alpha\mid B,b)\mathbb{I}_{[Q_j=1]}}{\mbox{Pr}_{\alpha}}$$

\noindent This makes the contrast estimators:\\
$\widehat{DE}(\alpha\mid B,b)= \widehat{\overline{Y}}(1;\alpha\mid B,b) - \widehat{\overline{Y}}(0;\alpha\mid B,b)$,\\
$\widehat{OE}(\alpha, \gamma\mid B,b)= \widehat{\overline{Y}}(\alpha\mid B,b) - \widehat{\overline{Y}}(\gamma\mid B,b)$,\\
$\widehat{IE}(\alpha, \gamma\mid B,b)= \widehat{\overline{Y}}(0;\alpha\mid B,b) - \widehat{\overline{Y}}(0;\gamma\mid B,b)$, and\\
$\widehat{TE}(\alpha, \gamma\mid B,b)= \widehat{\overline{Y}}(1;\alpha\mid B,b) - \widehat{\overline{Y}}(0;\gamma\mid B,b)$.\\

Let the individual level baseline variables conditional outcome estimators be defined by: 
$$\widehat{\overline{Y}}_j(1;\alpha\mid b) = \frac{1}{\sum_{i=1} \mathbb{I}_{[w_{ij} \in b]}}\sum^{n_j}_{i} \frac{Y_{ij} \mathbb{I}_{[z_{ij}=1]} \mathbb{I}_{[w_{ij} \in b]}}{\mbox{P}^{j}_{\alpha}}$$
with, $\widehat{\overline{Y}}_j(0;\alpha\mid b)$ similarly defined, and 
$$\widehat{\overline{Y}}_j(\alpha\mid b) = \frac{1}{\sum_{i=1} \mathbb{I}_{[w_{ij} \in b]}}\sum^{n_j}_{i} Y_{ij} \mathbb{I}_{[w_{ij} \in b]}.$$ The the population level estimators are given by 
$$\widehat{\overline{Y}}(z; \alpha\mid b) = \frac{1}{\sum_j \mathbb{I}_{[\sum_i \mathbb{I}_{[w_{ij}\in b]}>0]}} \frac{\sum\widehat{\overline{Y}}_j(z; \alpha\mid b) \mathbb{I}_{[Q_j=1]}}{\mbox{Pr}_{\alpha}},$$
and
$$\widehat{\overline{Y}}(\alpha\mid b) = \frac{1}{\sum_j \mathbb{I}_{[\sum_i \mathbb{I}_{[w_{ij}\in b]}>0]}}\frac{\sum\widehat{\overline{Y}}_j( \alpha\mid b) \mathbb{I}_{[Q_j=1]}}{\mbox{Pr}_{\alpha}},$$
with the contrast estimator following in the same way as above. 

Let the cluster level baseline variables conditional outcome estimators be defined by: 
$$\widehat{\overline{Y}}_j(1; \alpha\mid B) = \frac{1}{n_{j}}\frac{\sum^{n_j}_{i} Y_{ij}\mathbb{I}_{[Z_{ij}=1]}\mathbb{I}_{[d_j\in B]}}{\mbox{P}^{j}_{\alpha}}$$ with, $\widehat{\overline{Y}}_j(0; \alpha\mid B)$ similarly defined, and $$\widehat{\overline{Y}}_j(\alpha\mid B) = \frac{1}{n_j}\sum^{n_j}_{i} Y_{ij}\mathbb{I}_{[d_j\in B]}.$$ The the population level estimators are given by 
$$\widehat{\overline{Y}}(z;\alpha\mid B) = \frac{1}{\sum^{J}_{j=1}\mathbb{I}_{[d_j \in B]}}\sum^{J}_{j=1}\frac{\widehat{\overline{Y}}_j(z; \alpha\mid B)\mathbb{I}_{[Q_j=1]}}{\mbox{Pr}_{\alpha}}$$ 
and 
$$\widehat{\overline{Y}}(\alpha\mid B) = \frac{1}{\sum^{J}_{j=1}\mathbb{I}_{[d_j \in B]}}\sum^{J}_{j=1}\frac{\widehat{\overline{Y}}_j( \alpha\mid B)\mathbb{I}_{[Q_j=1]}}{\mbox{Pr}_{\alpha}}.$$ 
Again, the contrast estimator following in the same way as above. All of the $\gamma$ estimators are similarly defined. 

\textbf{RESULTS}
\begin{enumerate}
\item Under Assumptions a-c and when $M_{j,b}>0$
\begin{itemize}
    \item $\mbox{E}\{[\widehat{\overline{Y}}_j(z; \alpha\mid b)\mid Q_j=1]\} = \overline{Y}_j(z; \alpha\mid b)$
    \item $\mbox{E}\{[\widehat{\overline{Y}}_j(\alpha\mid b)\mid Q_j=1]\} = \overline{Y}_j(\alpha\mid b)$ 
    \item $\mbox{E}\{\widehat{DE}_j(\alpha\mid b)\mid Q_j=1\} = DE_j(\alpha\mid b)$.
\end{itemize}

\item Under assumptions a-c and when $M_{b}>0$
\begin{itemize}
    \item $\mbox{E}\{[\widehat{\overline{Y}}(\alpha\mid b)]\} = \overline{Y}(\alpha\mid b)$

\item $\mbox{E}\{[\widehat{\overline{Y}}(z; \alpha\mid b)]\} = \overline{Y}(z; \alpha\mid b)$

\item $\mbox{E}\{\widehat{DE}(\alpha\mid b)\} = DE(\alpha\mid b)$

\item $\mbox{E} \left\{\widehat{OE}(\alpha, \gamma\mid b)\right\} = OE(\alpha, \gamma\mid b)$

\item $\mbox{E} \left\{\widehat{IE}(\alpha, \gamma\mid b)\right\} = IE(\alpha, \gamma\mid b)$

\item $\mbox{E} \left\{\widehat{TE}(\alpha, \gamma\mid b)\right\} = TE(\alpha, \gamma\mid b)$.
\end{itemize}
\item Under Assumptions a-c
\begin{itemize}
    \item $\mbox{E}\{[\widehat{\overline{Y}}_j(z;\alpha\mid B)\mid Q_j=1]\} = \overline{Y}_j(z;\alpha\mid B)$ 

    \item$\mbox{E}\{[\widehat{\overline{Y}}_j(\alpha\mid B)\mid Q_j=1]\} = \overline{Y}_j(\alpha\mid B)$

    \item$\mbox{E}\{\widehat{DE}_j(\alpha\mid B)\mid Q_j=1\} = DE_j(\alpha\mid B)$
\end{itemize}
\item Under assumptions a-c and when $M_{B}>0$
\begin{itemize}
    \item $\mbox{E} \left\{\widehat{\overline{Y}}(z;\alpha\mid B)\right\} = \overline{Y}(z;\alpha\mid B)$ and
\item $\mbox{E} \left\{\widehat{\overline{Y}}(\alpha\mid B) \right\} = \overline{Y}(\alpha\mid B)$

\item $\mbox{E} \left\{\widehat{DE}(\alpha\mid B)\right\} = DE(\alpha\mid B)$

\item $\mbox{E} \left\{\widehat{OE}(\alpha, \gamma\mid B)\right\} = OE(\alpha, \gamma\mid B)$

\item $\mbox{E} \left\{\widehat{IE}(\alpha, \gamma\mid B)\right\} = IE(\alpha, \gamma\mid B)$

\item $\mbox{E} \left\{\widehat{TE}(\alpha, \gamma\mid B)\right\} = TE(\alpha, \gamma\mid B)$.
\end{itemize}
\item Under assumptions a-c and when $M_{j,b}>0$ 
\begin{itemize}
    \item $\mbox{E} \left\{[\widehat{\overline{Y}}_j(z;\alpha\mid B,b)\mid Q_j=1]\right\} =\overline{Y}_j(z;\alpha\mid B,b)$

\item$\mbox{E} \left\{[\widehat{\overline{Y}}_j(\alpha\mid B,b)\mid Q_j=1]\right\} =\overline{Y}_j(\alpha\mid B,b)$

\item$\mbox{E}\{\widehat{DE}_j(\alpha\mid B,b)\mid Q_j=1\} = DE_j(\alpha\mid B,b)$.
\end{itemize}
\item Under assumptions a-c and when $M_{B,b}>0$
\begin{itemize}
    \item $\mbox{E} \left\{\widehat{\overline{Y}}(z;\alpha\mid B,b)\right\} = \overline{Y}(z;\alpha\mid B,b)$ 
 \item$\mbox{E} \left\{\widehat{\overline{Y}}(\alpha\mid B,b) \right\} = \overline{Y}(\alpha\mid B,b)$

  \item $\mbox{E} \left\{\widehat{DE}(\alpha\mid B,b)\right\} = DE(\alpha\mid B,b)$

  \item $\mbox{E} \left\{\widehat{OE}(\alpha, \gamma\mid B,b)\right\} = OE(\alpha, \gamma\mid B,b)$

  \item $\mbox{E} \left\{\widehat{IE}(\alpha, \gamma\mid B,b)\right\} = IE(\alpha, \gamma\mid B,b)$

  \item $\mbox{E} \left\{\widehat{TE}(\alpha, \gamma\mid B,b)\right\} = TE(\alpha, \gamma\mid B,b)$.
  \end{itemize}
\end{enumerate}

\textbf{Theorem 1}\\ 
Under assumptions a-d and $M_{j,b}>0$ and $n_j P_{\alpha}-1>0$, 
$\mbox{E}[\widehat{\mbox{Var}} 
(\widehat{\overline{Y}}_j(z;\alpha\mid B,b)\mid Q_j=1)\mid Q_j=1] = \mbox{Var}(\widehat{\overline{Y}}_j(z;\alpha\mid b)\mid Q_j=1),$
where 
\begin{eqnarray*}
&&\widehat{\mbox{Var}}(\widehat{\overline{Y}}_j(1;\alpha\mid B,b)\mid Q_j=1) \equiv\\
&&\mathbb{I}_{[d_j \in B]}\left(1-P^{j}_{\alpha}\right)\frac{\sum^{n_j}_{i=1} Z_{ij}\{Y_{ij}Z_{ij}\mathbb{I}_{[w_{ij} \in b]}(n_j/M_{j,b}) - \widehat{\overline{Y}}_j(1;\alpha\mid b)\}^2}{(n_j P^{j}_{\alpha}-1)n_j P^{j}_{\alpha}}
\end{eqnarray*}
with $\widehat{\mbox{Var}}(\widehat{\overline{Y}}_j(0;\alpha\mid B,b)\mid Q_j=1)$ and 
$\widehat{\mbox{Var}}(\widehat{\overline{Y}}_j(z;\gamma\mid B,b)\mid Q_j=0)$ defined similarly.

\textbf{Theorem 2}\\
Under assumptions a-d and $M_{B,b}>0$ $n Pr_{\alpha}-1>0$ and $n_j P_{\alpha}-1>0$ for all $j$, 

$\mbox{E}[\widehat{\mbox{Var}} 
(\widehat{\overline{Y}}(z;\alpha\mid B,b))] = \mbox{Var}(\widehat{\overline{Y}}(z;\alpha\mid B,b)),$
where 
\begin{eqnarray*}
\widehat{\mbox{Var}}(\widehat{\overline{Y}}(z;\alpha\mid B,b)) 
&\equiv&
\left(1- Pr_{\alpha}\right)\frac{\sum^{J}_{j=1} \mathbb{I}_{[Q_{j}=1]}\{\widehat{\overline{Y}}_{j}(z;\alpha\mid B,b)(J/M_{B,b}) - \widehat{\overline{Y}}(z;\alpha\mid B,b)\}^2}{(J Pr_{\alpha}-1) J Pr_{\alpha}}\\
&+&
\frac{1}{Pr_{\alpha} M^2_{B,b}}\sum_{j \in J^b}\widehat{\mbox{Var}}[\widehat{Y}_{j}(z;\alpha\mid B,b)]\mathbb{I}_{[Q_{j}=1]}
\end{eqnarray*}
with $\widehat{\mbox{Var}}(\widehat{\overline{Y}}(z;\gamma\mid B,b))$ defined similarly. Proofs of all theorems and results are given in the supplementary materials. The proofs of theorems 1 and 2 consider conditioning under each type, group and individual level subgroups, separately as well as together. 

\section{Numerical Example} \label{Example}

\begin{table}[h!]
\caption{Example Data\label{table1}}
\begin{tabular}{cccccc|cccccc}
group&ID&Y(1)&Y(0)&$w_{ij}\in$b&$d_j\in$B&group&ID&Y(1)&Y(0)&$w_{ij}\in$b&$d_j\in$B\\
1&11&3&0&1&1&2&21&0&2&0&0\\
1&12&2&0&0&1&2&22&2&3&0&0\\
1&13&10&2&1&1&2&23&4&6&0&0\\
1&14&1&1&0&1&2&24&5&7&0&0\\
\hline \\
3&31&1&2&1&0&4&41&0&3&1&0\\
3&32&2&1&0&0&4&42&2&1&0&0\\
3&33&3&0&1&0&4&43&4&5&0&0\\
3&34&10&1&0&0&4&44&5&7&0&0\\
\end{tabular}
\footnote{Subjects only have two counterfactual outcomes that need to be considered because of our assumption of stratified interference.}
\end{table}

Consider a setting in which there are four groups of four subjects each, in which two of four groups will received 50\% coverage ($\alpha=0.5$, exactly 2 of 4 individuals per group receive treatment), and the other two of four groups will receive 25\% coverage ($\gamma=0.25$, exactly 1 of 4 individuals per group receive treatment). Example data are given in Table \ref{table1}. There are 6 total ways to randomize the clusters to 50\% or 25\%  in a one-to-one ratio. Let the realized randomization at the group level be, $q=\{0,1,0,1\}$ and that I wish to do a subgroup analysis based on the individual level covariate $b$. For example, this group may represent the sex of the subjects, which may not normally be rare, but could be in some settings. Under this randomization, there are six possible randomizations for group 4, $\mathbf{z}_4$ with 50\% coverage, $\{0,1,1,0\}$, $\{0,0,1,1\}$, $\{0,1,0,1\}$, $\{1,1,0,0\}$, $\{1,0,0,1\}$, $\{1,0,1,0\}$. By the definition of the true group average value given above, $\overline{Y}_4(0; \alpha\mid b) = 3$, as subject 41 has a value under placebo of 3. This also makes clear that the estimand definitions in this setting are sensible, as you would want to know the average value within the subgroup had all subjects within the group been assigned to placebo, but still under 50\% coverage.

I want to estimate $\overline{Y}_4(0; \alpha\mid b)$, the individual variable conditional group average outcome for the untreated. Let us consider the possible estimators. One possible ``natural" estimator for the individual level conditional group average outcome under $z=1$ based on \citet{hudgens08} would be 
\begin{equation}
 \widehat{\overline{Y}}^{N}_j(0; \alpha\mid b)  = \frac{1}{\sum (1-z_{ij}) \mathbb{I}_{[w_{ij} \in b]}} \sum Y_{ij} (1-z_{ij}) \mathbb{I}_{[w_{ij} \in b]}, 
\end{equation}
where the superscript $N$ is for natural. 

Similarly, I consider a possible subgroup conditional version of the H\'ajek estimator \citet{hakek71}, 
\begin{equation}
 \widehat{\overline{Y}}^{Hj}_j(0; \alpha\mid b)  =  \frac{\sum \frac{\mathbb{I}_{[z_{ij}=0]} \mathbb{I}_{[w_{ij} \in b]}}{1-P^{j}_{\alpha}} Y_{ij}}{\sum \frac{\mathbb{I}_{[z_{ij}=0]} \mathbb{I}_{[w_{ij} \in b]}}{1-P^{j}_{\alpha}}}.
\end{equation}

Finally, the estimator given above is: 
\begin{equation}
\widehat{\overline{Y}}_j(0;\alpha\mid b) = \frac{1}{\sum_i \mathbb{I}_{[d_{j} \in b]}}\sum^{n_j}_{i} \frac{Y_{ij}\mathbb{I}_{[Z_{ij}=0]}\mathbb{I}_{[w_{ij} \in b]}}{1-\mbox{P}^{j}_{\alpha}}.
\end{equation}

In group 4, the  ``natural" and  H\'ajek estimators will either be undefined, denoted as ``NA" in the table, or have a different but defined value than the HT estimator if the one subject with $w_{ij} \in b$ is assigned to treatment, as is displayed in Table \ref{table2}. If I instead set the  ``natural" and  H\'ajek estimators to zero when they were undefined, this would result in a bias towards 0. As defined in theorem 1, $\widehat{\mbox{Var}}(\widehat{\overline{Y}}_4(0;\alpha\mid b)\mid Q_j=1)$ will be equal to, 18, 0, 0, 18, 0, 18, giving an average of 9, which is the true variance of $ \widehat{\overline{Y}}_4(0;\alpha\mid b)$ over the possible randomizations, as can be seen in Table \ref{table2}. 

\begin{table}[h!]
\caption{{Estimator Values over the possible randomizations of Group 4}\label{table2}}
\begin{tabular}{l|cccccc|l}
&\multicolumn{6}{c}{Randomizations}\\
 &$\{0,1,1,0\}$&$\{0,0,1,1\}$& $\{0,1,0,1\}$& $\{1,1,0,0\}$& $\{1,0,0,1\}$& $\{1,0,1,0\}$& Average\\
 \hline
$ \widehat{\overline{Y}}^{N}_4(0; \alpha\mid b)$&3&NA&NA&3&NA&3&NA\\
 \hline
 $ \widehat{\overline{Y}}^{Hj}_4(0; \alpha\mid b)$&3&NA&NA&3&NA&3&NA\\
  \hline
 $ \widehat{\overline{Y}}_4(0;\alpha\mid b)$&6&0&0&6&0&6&3\\
  \hline
\end{tabular}
\end{table}

If instead one wanted to estimate $\overline{Y}_3(0; \gamma \mid b)$, the individual variable conditional group average outcome for the untreated under 25\% coverage, all estimators would be defined and unbiased. Note that, the estimators are the same as given above for group 4 but with $\mbox{P}^{j}_{\alpha}$ replaced with $\mbox{P}^{j}_{\gamma}$. Let us consider the variance of the HT estimator. As defined in theorem 1, $\widehat{\mbox{Var}}(\widehat{\overline{Y}}_3(0;\gamma \mid b)\mid Q_j=0)$ will take on the values, 0, 0.44, 0.44, 0.44, over the randomizations, for an average of 0.33, which can be seen to be the true variance of $\widehat{\overline{Y}}_3(0;\gamma \mid b)$ over the 4 randomizations in Table \ref{table3}. This is a lower variance than either of the other two estimators, which both have a true variance over the randomizations of 0.5. 

\begin{table}[h!]
\caption{{Estimator Values over the possible randomizations of Group 3}\label{table3}}
\begin{tabular}{l|cccc|l}
 &$\{1,0,0,0\}$&$\{0,1,0,0\}$& $\{0,0,1,0\}$& $\{0,0,0,1\}$&Average\\
  \hline
$ \widehat{\overline{Y}}^{N}_3(0; \gamma \mid b)$&0&1&2&1&1\\
 \hline
 $ \widehat{\overline{Y}}^{Hj}_3(0; \gamma\mid b)$&0&1&2&1&1\\
  \hline
 $ \widehat{\overline{Y}}_3(0;\gamma\mid b)$&0&1.333&1.333&1.333&1\\
  \hline
\end{tabular}
\end{table}

I do not further speculate on how the other estimators would be extended to the cross-conditioning setting. However, in this setting the HT style estimator $ \widehat{\overline{Y}}(1; \alpha\mid B, b)$, for example, under the realized randomization $q=\{0,1,0,1\}$, would be 0 and the variance estimate would be 0, as group 1, the only group with $d_j \in B$, is assigned to $\gamma$ coverage.

\section{Discussion} \label{DIS}
In this article I define a set of conditional estimators for rare subgroup analysis following the style of the HT marginal estimator. I show that these estimators are unbiased provided there is at least one subject in the group, or population, of interest within the subgroup of interest. I provide variance estimates for average groups level and population level subgroup conditional outcomes under treatment or control for a given coverage level, and prove that they are unbiased over the possible randomizations. 

In this paper I point out a useful characteristic of Horvitz-Thompson style estimators that the denominator of the estimator is not dependent on the randomization and therefore these estimators can be pre-specified for use in subgroups analysis in blinded clinical trials without a fear that they will be undefined due to randomization. This does not hold for other estimators that have been suggested for estimation in the presence of interference, although there maybe ways to modify them to allow for subgroup analysis. In larger or less rare subgroups where, for the given randomization strategy, it is not possible to have zero subjects from the subgroup within a randomized arm, both the ``natural" estimator and the H\'ajek estimators will be unbiased. 

I make a large number of assumptions, and has been pointed out in the long literature since \citet{hudgens08}, estimation and some inference is still possible under weaker assumptions \citet{aronow2017estimating, savje2017average}. It is a future area of research to consider subgroup analysis under fewer or more relaxed assumptions.

\bibliographystyle{plainnat}
\bibliography{clusterlevel}

\begin{thebibliography}{14}
\providecommand{\natexlab}[1]{#1}
\providecommand{\url}[1]{\texttt{#1}}
\expandafter\ifx\csname urlstyle\endcsname\relax
  \providecommand{\doi}[1]{doi: #1}\else
  \providecommand{\doi}{doi: \begingroup \urlstyle{rm}\Url}\fi

\bibitem[Aronow et~al.(2017)Aronow, Samii, et~al.]{aronow2017estimating}
Peter~M Aronow, Cyrus Samii, et~al.
\newblock Estimating average causal effects under general interference, with
  application to a social network experiment.
\newblock \emph{The Annals of Applied Statistics}, 11\penalty0 (4):\penalty0
  1912--1947, 2017.

\bibitem[H\'ajek(1971)]{hakek71}
J.~H\'ajek.
\newblock Comment on a paper by {D. Basu} in {F}oundations of {S}tatistical
  {I}nference , {G}odambe {V}. \& {S}prott {D}. eds., 1971.

\bibitem[Halloran and Hudgens(2012)]{halloran12}
M~Elizabeth Halloran and Michael~G Hudgens.
\newblock Causal inference for vaccine effects on infectiousness.
\newblock \emph{The International Journal of Biostatistics}, 8\penalty0
  (2):\penalty0 1--40, 2012.

\bibitem[Horvitz and Thompson(1952)]{Horvitz52}
D.~Horvitz and D.~Thompson.
\newblock A generalization of sampling without replacement from a finite
  universe.
\newblock \emph{Journal of the American Statistical Association}, 47\penalty0
  (260):\penalty0 663--685, 1952.

\bibitem[Hudgens and Halloran(2008)]{hudgens08}
Michael~G Hudgens and M~Elizabeth Halloran.
\newblock Toward causal inference with interference.
\newblock \emph{Journal of the American Statistical Association}, 103\penalty0
  (482):\penalty0 832--842, 2008.

\bibitem[Liu and Hudgens(2014)]{liu14}
Lan Liu and Michael~G Hudgens.
\newblock Large sample randomization inference of causal effects in the
  presence of interference.
\newblock \emph{Journal of the American Statistical Association}, 109\penalty0
  (505):\penalty0 288--301, 2014.

\bibitem[Liu et~al.(2016)Liu, Hudgens, and Becker-Dreps]{liu2016inverse}
Lan Liu, Michael~G Hudgens, and Sylvia Becker-Dreps.
\newblock On inverse probability-weighted estimators in the presence of
  interference.
\newblock \emph{Biometrika}, 103\penalty0 (4):\penalty0 829--842, 2016.

\bibitem[S{\"a}vje et~al.(2017)S{\"a}vje, Aronow, and
  Hudgens]{savje2017average}
Fredrik S{\"a}vje, Peter~M Aronow, and Michael~G Hudgens.
\newblock Average treatment effects in the presence of unknown interference.
\newblock \emph{arXiv preprint arXiv:1711.06399}, 2017.

\bibitem[Sobel(2006)]{sobel06}
Michael~E Sobel.
\newblock What do randomized studies of housing mobility demonstrate? causal
  inference in the face of interference.
\newblock \emph{Journal of the American Statistical Association}, 101\penalty0
  (476):\penalty0 1398--1407, 2006.

\bibitem[Tchetgen and VanderWeele(2012)]{tchetgen10}
Eric J~Tchetgen Tchetgen and Tyler~J VanderWeele.
\newblock On causal inference in the presence of interference.
\newblock \emph{Statistical methods in medical research}, 21\penalty0
  (1):\penalty0 55--75, 2012.

\bibitem[VanderWeele(2009)]{vanderweele09}
Tyler~J VanderWeele.
\newblock Concerning the consistency assumption in causal inference.
\newblock \emph{Epidemiology}, 20\penalty0 (6):\penalty0 880--883, 2009.

\bibitem[VanderWeele(2012)]{vanderweele12}
Tyler~J VanderWeele.
\newblock Mediation analysis with multiple versions of the mediator.
\newblock \emph{Epidemiology (Cambridge, Mass.)}, 23\penalty0 (3):\penalty0
  454--463, 2012.

\bibitem[VanderWeele(2013)]{vanderweele13}
Tyler~J VanderWeele.
\newblock Mediation and spillover effects in group- randomized trials: A case
  study of the 4rs educational intervention.
\newblock \emph{Biometrics}, 69\penalty0 (3):\penalty0 561--565, 2013.

\bibitem[VanderWeele and Tchetgen(2011)]{vanderweele11}
Tyler~J VanderWeele and Eric J~Tchetgen Tchetgen.
\newblock Effect partitioning under interference in two-stage randomized
  vaccine trials.
\newblock \emph{Statistics \& Probability Letters}, 81\penalty0 (7):\penalty0
  861--869, 2011.

\end{thebibliography}

\clearpage

\end{document}